\documentclass[11pt, reqno]{amsart}

\usepackage{color}

 \usepackage{amsmath}
 \usepackage{amssymb}
\usepackage{bm}

\newcommand{\mysection}[1]{\section{#1}
      \setcounter{equation}{0}}

\newtheorem{theorem}{Theorem}[section]

\newtheorem{corollary}[theorem]{Corollary}

\theoremstyle{definition}

\theoremstyle{remark}
\newtheorem{remark}{Remark}[section]

\newcommand{\tr}{\text{\rm tr}\,}

\newcommand\sfu{{\sf u}}

 \makeatletter
 \def\dashint{%
 \operatorname%
 {\,\,\text{\bf--}\kern-.98em\DOTSI\intop\ilimits@\!\!}}
\def\ninf{\qopname\relax\@empty{inf\phantom{p}\!\!\!}}
\def\dashnorm{\,\,\text{\bf--}\kern-.5em\|}
 \makeatother

\newcommand\bbeta{\text{\raise-.2ex\hbox{$\bm{\beta}$}}}

\newcommand\bR{\mathbb{R}}

\newcommand\bS{\mathbb{S}}

\begin{document}

\title[Elliptic equations with $L_{d}$-drift]
{Aleksandrov's estimates for elliptic equations with
drift in Morrey spaces containing $L_{d}$}

\author{Hongjie Dong}
\address[H. Dong]{Division of Applied Mathematics, Brown University,
182 George Street, Providence, RI 02912, USA}
\thanks{H. Dong is partially supported by the Simons Foundation, grant \# 709545.}
\email{Hongjie\_Dong@brown.edu}

\author{N.V. Krylov}
\address[N. V. Krylov]{127 Vincent Hall, University of Minnesota,
 Minneapolis, MN, 55455}
\email{nkrylov@umn.edu}



\begin{abstract}
 In this note, we obtain a version of Aleksandrov's maximum principle when the drift coefficients are in Morrey spaces, which contains $L_d$, and when the free term is in $L_p$ for some $p<d$.
\end{abstract}

\maketitle

\mysection{Introduction and main result}

Let $\bR^{d}$ be a Euclidean space of points
$x=(x^{1},...,x^{d})$. Define
$$
D_{i}=\partial/\partial x^{i},\quad D=(D_{i}),
\quad D_{ij}=D_{i}D_{j},\quad D^{2}=(D_{ij}),
$$
$$
B_{r}(x)=\big\{y\in\bR^{d}:|x-y|<r\big\},\quad B_{r}=B_{r}(0).
$$

Let $\bS$ be the set of symmetric $d\times d$
matrices and, for
$\delta\in(0,1]$,  let
$$
\bS_{\delta}=\big\{a\in\bS:\delta|\xi|^{2}\leq a^{ij}\xi^{i}\xi^{j}
\leq\delta^{-1}|\xi|^{2},\quad \forall \xi\in\bR^{d}\big\}.
$$

We take a measurable $\bS_{\delta}$-valued function $a$
on $B_{1}$ and a measurable $\bR^{d}$-valued function $b$
on $B_{1}$ and set
$$
L=a^{ij}D_{ij}+b^{i}D_{i}.
$$

The goal of this article is to prove the
following result in which $d_{0}=d_{0}(d,\delta)\in(d/2,d)$
is specified later.

\begin{theorem}
                                       \label{theorem 2.26.1}
Let $R\in(0,\infty)$, $d_{0}<q<d$, $p\in (d_{0}d/q,d)$,
$v\in W^{2}_{p}(B_{R})$ and
assume that for any $x\in B_{R}$ and $\rho\leq2R$
\begin{equation}
                                 \label{2.26.3}
\frac{1}{|B_{\rho}(x)|}\int_{B_{\rho}(x)}
I_{B_{R}}(y)|b|^{q}(y)\,dy\leq \bar b \rho^{-q},
\end{equation}
where $\bar b>0$ is sufficiently small
depending only on $d,\delta,p,q$ (independent of $R$).
Then in $B_{R}$ we have
\begin{equation}
                                 \label{2.26.2}
v\leq \sup_{\partial B_{R}}v
+NR^{2-d/p}\|(Lv)_{-}\|_{L_{p}(B_{R})},
\end{equation}
where    $N$
depends only on $d,\delta,p$, and $q$.

\end{theorem}

The first results about estimates of the maximum
of
solutions of elliptic equations
even not uniformly elliptic in terms of
the $L_{d}$-norms of the free terms belong to
A.D. Aleksandrov in 1960, see Theorem 8 in
\cite{Al_60}, where he already
considered $b\in L_{d}$.
The proofs are given in 1963 in \cite{Al_63}.

There was considerable interest in reducing $L_{d}$-norm
of the free term to $L_{d_{0}}$ norm with
$d_{0}<d$. This was achieved by Cabr\'e
\cite{Ca_95} for bounded $b$
and by Fok in \cite{Fo_98} for
  $b\in L_{d+\varepsilon}$.
In
 \cite{Kr_21}   the second named author
allowed $b\in L_{d}$ and the free term  in
$L_{d_{0}}$   with
$d_{0}<d$. This made it possible
to develop  in  \cite{Kr_21b}   a $W^{2}_{p}$-solvability
theory for linear equations with $b\in L_{d}$ and $p<d$
(without assuming that the $L_{d}$-norm of $b$ is
small or the domain is small). Applied
to fully nonlinear equations we can now
treat $W^{2}_{d_{0}}$-solvability with
the ``coefficients'' of the first order terms
in $L_{d}$ (see \cite{Kr_20}).
 In this paper we are doing
the next step by reducing the integrability
property of $b$ further down.

\begin{remark}
                                     \label{remrk 2.26.1}
Observe that \eqref{2.26.3} may be satisfied but
$b\not\in L_{d}$. Our arguments are based on the theory of PDEs. By using probabilistic means, it is proved in \cite{Kr_21a} that once
condition \eqref{2.26.3} is satisfied in all balls
 with  radius  $\rho\leq   \rho_0 $  for some small $\rho_0$,
 estimate
\eqref{2.26.2} holds for all $R$ with $N$
that also depends on $R$  and $\rho_0$. For $R$ small,
condition \eqref{2.26.3} is indeed
satisfied if $b\in L_{d}$ because by H\"older's
inequality the left-hand side of \eqref{2.26.3}
is less than $N(d)\|b\|{ ^q}_{L_{d}(B_{R})} \rho^{-q }$. It is then
tempting to claim that our Theorem \ref{theorem 2.26.1}
contains Aleksandrov's result specified for
uniformly nondegenerate equations with bounded $a$
or contains the corresponding result of
\cite{Kr_21}
(see Corollary 3.1 there).
However in Aleksandrov's result and in \cite{Kr_21}  the constants depend
only on the $L_{d}$-norm of $b$  and not on the
rate with which $\|b\|_{L_{d}(B_{R}(x))}\to0$
as $R\downarrow 0$.

On the other hand, one cannot drop the smallness
assumption of $\bar b$. For instance,
if $R=1$, $a^{ij}=\delta^{ij}$,
 $b(x)=-cx/|x|^{2}$ and $v(x)=1-|x|^{2}$,
then $Lu=-2d+2c$ and if $c\geq d$, $Lu\geq0$
and \eqref{2.26.2} is false.

It turns out that
\eqref{2.26.2} is true
for any $u\in W^{2}_{p}(B_{1})$ if $c$ is small enough.
Indeed,
observe that if $|x|\geq 2\rho$ the integral
$$
 \int_{B_{\rho}(x)}
I_{B_{1}}(y)|y|^{-q}\,dy
$$
is less than $\rho^{-q}|B_{\rho}(x)|$
and if $|x|\leq 2\rho$, the said integral is dominated by
$$
\int_{B_{3\rho} }
 |y|^{-q}\,dy=N(d)\rho^{n-q}.
$$

This example
has a deep connection with the so-called
form-boundedness condition
from \cite{KS_20}. We also stress one more time
that in this example $b\not\in L_{d}$.

\end{remark}

\mysection{Proof of Theorem \protect\ref{theorem 2.26.1}}

By using change of scale we see that it suffices
to concentrate on $R=1$. However since we also need
$B_{2}$, we keep $R$ free for a while.

For $p\in[1,\infty)$ and $\mu\in(0,d/p)$, introduce
Morrey's space $E_{p,\mu}(B_{R})$
as the set of $g\in L_{p}(B_{R})$ such that
$$
\|g\|_{E_{p,\mu}(B_{R})}:=
\sup_{\rho\leq 2R,x\in B_{R}}\rho^{\mu}
\dashnorm g \|_{L_{p}(B_{R,\rho}(x))}<\infty,
$$
where $B_{R,\rho}(x)=B_{R}\cap B_{\rho}(x)$
and
$$
\dashnorm g\|_{L_{p}(\Gamma)}=\Big(\frac{1}{|\Gamma|}
\int_{\Gamma}|g|^{p}\,dx\Big)^{1/p}.
$$
Let
$$
E^{2}_{p,\mu}(B_{R})=\{u:u,Du,D^{2}u\in E_{p,\mu}(B_{R})\}.
$$
and provide $E^{2}_{p,\mu}(B_{R})$ with an obvious norm.
The case $R=\infty$ is not excluded and we drop
$B_{\infty}=\bR^{d}$ from our notation.

\begin{theorem}
                                       \label{theorem 2.24.1}
Let $1<\mu\leq d/p$, $p>1$. Define $q$ ($>p$) from
$$
\frac{1}{q}=\frac{1}{p}-\frac{1}{\mu p}.
$$
Then for any $u\in W^{1}_{p}(\bR^{d})$,
\begin{equation}
                                               \label{2.24.1}
\|u\|_{E_{q,\mu p/q}}\leq N(d,p,\mu)\|Du\|_{E_{p,\mu}}.
\end{equation}
\end{theorem}

Proof. As it follows from Secs.~1, 2, Ch.~V of \cite{St_70},
for almost any $x$ we have
$$
|u(x)|\leq N(d)\int_{\bR^{d}}|Du(y)||x-y|^{-d+1}\,dy.
$$
After that \eqref{2.24.1} follows from
Theorem 3.1 of \cite{Ad_75}. The theorem is proved.

This theorem is quite remarkable because it allows us
to estimate higher powers of $u$ compared with the usual
Sobolev embedding theorem at the expense of requiring
$Du$ be slightly better. It is the first crucial point
in the proof of Theorem \ref{theorem 2.26.1}.
Another one, well expected, is Theorem \ref{theorem 2.24.2}.

\begin{corollary}
                                     \label{corollary 2.24.1}
Let $1<p<q< d$  and $b\in E_{q,1}(B_{ 1})$.
Set $\mu=q/p$.
Then for any $u\in E^{2}_{p,\mu}(B_{2})$, we have
$$
\|b|D u|\|_{E_{p,\mu}(B_{1})}
\leq N\|b\|_{E_{q,1}(B_{1})}
 \|u\|_{E^{2}_{p,\mu}(B_{2 })},
$$
where the constants $N$ depend only on $d,p,q$.
\end{corollary}

Proof. Take $x\in B_{1}$, $\rho\leq 2$,
and take $\zeta\in C^{\infty}_{0}(\bR^{d})$
such that $\zeta=1$ on $B_{1}$, $\zeta=0$
outside $B_{2 }$, and $|\zeta|+|D\zeta|+|D^{2}\zeta|
\leq N=N(d)$.
By using H\"older's inequality, we see that
$$
\rho^{\mu}
\dashnorm b|Du| \|_{L_{p}(B_{1,\rho}(x))}
\leq N\rho\dashnorm b \|_{L_{q}(B_{1,\rho}(x))}
\rho^{\mu-1}\dashnorm D(\zeta u)\|_{L_{q'}(B_{\rho}(x))},
$$
where $q'=pq/(q-p)$ and the constant $N$ arose
because $|B_{1,\rho}(x)|$ is not quite $|B_{\rho}(x)|$.
Furthermore, since $\mu-1=\mu p/q'$ and $1/q'=1/p-1/(\mu p)$,
by Theorem \ref{theorem 2.24.1}
$$
\rho^{\mu-1}\dashnorm D(\zeta u)\|_{L_{q'}( B_{ \rho}(x))}
\leq N\|D^{2}(\zeta u)\|_{E_{p,\mu}}
\leq N\|u\|_{E^{2}_{p,\mu}(B_{2 })}.
$$
This obviously leads to the desired result.

For $\sfu''\in\bS$ introduce a Pucci function
$$
P(\sfu'')=\sup_{a\in\bS_{\delta}}\tr (a\sfu'').
$$

In the following, by $d_{0}$ we denote the constant called
$n_{0}$ in \cite{BLP_16} corresponding to the domain
$B_{2}$ and the operator $P(D^{2}u)$. Note
that $d_{0}=d_{0}(d,\delta)\in(d/2,d)$.

\begin{theorem}
                                     \label{theorem 2.24.2}
Let $d_{0}<p<q<d$ and set $\mu=q/p$. Assume that
a nonnegative {\em bounded\/}  function
$b\in E_{q,1}(B_{1})$ and $b=0$ outside $B_{1}$. Then there is
a $\bar b=\bar b(d,\delta,q,p)>0$ such that if
$\|b\|_{E_{q,1}(B_{1})}\leq \bar b$, then for any
$f\in E_{p,\mu}(B_{2 })$ there exists a unique
$u\in E^{2}_{p,\mu}(B_{2 })\cap C(\bar B_{2 })$ satisfying
\begin{equation}
                                      \label{2.24.4}
P(D^{2}u)+b|Du|+f=0
\end{equation}
in $B_{2}$ and equal  to  zero on $\partial B_{2}$.
Moreover, we have
\begin{equation}
                                      \label{2.24.40}
 \|u\|_{E^{2}_{p,\mu}(B_{2})}
\leq N\|f\|_{E_{p,\mu}(B_{2})},
\end{equation}
where $N$ depends only on  $d$, $\delta$, $q$, and $p$.

\end{theorem}

Proof. The existence and uniqueness of
solution follows directly from Theorem 4.1
of \cite{BLP_16} due to the assumption that $b$
is {\em bounded\/}. To prove \eqref{2.24.40},
it suffices to observe that by the same Theorem 4.1
of \cite{BLP_16}
$$
\|u\|_{E^{2}_{p,\mu}(B_{2})}\leq N\|b|Du|+f\|_{E_{p,\mu}(B_{2})}
\leq N\|b|Du| \|_{E_{p,\mu}(B_{2})}+N
\| f\|_{E_{p,\mu}(B_{2})},
$$
where, thanks to Corollary
\ref{corollary 2.24.1},
 the first term  on the right-hand side  can be absorbed into the left-hand side
if $\|b\|_{E_{p,1}(B_{1})}$ is sufficiently small.
The theorem is proved.

Now we prove Theorem \ref{theorem 2.26.1}
in case $R=1$ when it takes the following form.

\begin{theorem}
                                       \label{theorem 2.25.2}
Let $d_{0}<q<d$, $p\in (d_{0}d/q,d)$,  and
$v\in W^{2}_{p}(B_{1})$. Assume that
$b\in E_{q,1}(B_{1})$ and
$\|b\|_{E_{q,1}(B_{1})}\leq \bar b$,
where $\bar b$ is taken from Theorem \ref{theorem 2.24.2}.
Then in $B_{1}$ we have
\begin{equation}
                                 \label{2.26.1}
v\leq \sup_{\partial B_{1}}v
+N\|(Lv)_{-}\|_{L_{p}(B_{1})},
\end{equation}
 where $N$ depends only on $d,\delta,q$, and $p$.

\end{theorem}

Proof. If we replace $L$ in \eqref{2.26.1}   with
$I_{|b|>n}\Delta+I_{|b|\leq n}L$ and assume that
our assertion is true, then by passing to the limit
as $n\to\infty$ and using the dominated and monotone
convergence theorems, we obtain \eqref{2.26.1} as is.
It follows that we may assume that $b$ is bounded
and then that $v$, $b$, and $a$ are smooth.
After that by subtracting from $v$ the solution
$w$ of $Lw=0$ in $B_{1}$ with boundary value
$w=v$, we reduce the situation to the one
in which $v=0$ on $\partial B_{1}$. Then
set $f=(Lv)_{-}I_{B_{1}}$, extend $b$ as zero
outside $B_{1}$, and
define $u$ as a solution of
\eqref{2.24.4} in $B_{2}$ with zero boundary data.
According to \cite{Wi_09} or \cite{Kr_18}, there is such
  solution $u$ which belongs to
$  W^{2}_{r}(B_{2})$ for any $r>1$.
By the maximum principle, $u\geq0$ in $B_{2}$
and, since $L(u-v)\leq
Lu +f\leq0$ in $B_{1}$, again by the
maximum principle $v\leq u$ in $B_{1}$.

Furthermore, in light of H\"older's
inequality, $u\in E^{2}_{r,\nu}(B_{2 })$ for any $r \in (1,\infty) $ and $\nu\in (0,d/r)$.
Now, by embedding theorems, \eqref{2.24.40}, and again
H\"older's inequality, for $d_{0}<r<q$ and $\nu=q/r$,
$$
u\leq N\|u\|_{W^{2}_{r}(B_{2 })}
\leq N\|u\|_{E^{2}_{r,\nu}(B_{2 })}
\leq N\|f\|_{E_{r,\nu}(B_{2 })}\leq N
\|f\|_{L_{rd/q}(B_{2})}.
$$
After that it only remains to note that
$rd/q=p$ for $r=pq/d$ ($ r\in (d_{0},q)$ with $\nu=q/r=d/p>1$).
The theorem is proved.

\begin{remark}
As an intermediate result we have proved that if $v$, $a$, $b$
are smooth and $v=0$ on $\partial B_{1}$, then
$$
v\leq N\|(Lv)_{-}\|_{E_{qp/d,d/p}(B_{1})}.
$$
\end{remark}

{\bf Acknowledgment}.
The authors are sincerely grateful to Doyoon Kim
for pointing out several glitches in
the first version of the article.


\begin{thebibliography}{mm}

\bibitem{Ad_75} D. Adams, {\em A note on Riesz potentials\/},
Duke Math. J., Vol. 42 (1975), No. 4, 765-778.

\bibitem{Al_60} A. D. Aleksandrov, {\em
Certain estimates for the Dirichlet problem\/},
Dokl. Akad. Nauk SSSR, Vol. 134 (1960),
 1001--1004 (Russian); translated as
Soviet Math. Dokl., Vol. 1 (1961) 1151--1154.

\bibitem{Al_63} A. D. Aleksandrov,
{\em Uniqueness conditions and estimates for
the solution of the Dirichlet
problem\/}, Vestnik Leningrad. Univ., Vol. 18 (1963), No. 3,
5-29 in Russian;
English translation in
Amer. Mat. Soc. Transl., Vol. 68 (1968), No. 2, 89-119.


\bibitem{BLP_16} S.S. Byun, M. Lee, and D.K. Palagachev,
{\em Hessian estimates in weighted Lebesgue spaces for
 fully nonlinear
elliptic equations\/}, J. Differential Equations,
Vol. 260 (2016), No. 5, 4550-4571.



\bibitem{Ca_95}
X. Cabr\'e, {\em
On the Alexandroff-Bakelman-Pucci estimate and the reversed
H\"older inequality for solutions of elliptic
and parabolic equations\/},
Comm. Pure Appl. Math., Vol. 48 (1995), 539--570.

\bibitem{Fo_98} K. Fok, {\em A nonlinear Fabes-Stroock result\/},
Comm. PDEs, Vol 23 (1998), No. 5-6, 967--983.


 \bibitem{KS_20} D. Kinzebulatov and Yu.A. Semenov,
{\em Brownian motion with general drift\/},
Stoch. Proc. Appl., Vol. 130 (2020), No. 5, 2737-2750.

\bibitem{Kr_18} N.V. Krylov,
``Sobolev and viscosity solutions for fully nonlinear  elliptic
and parabolic equations'', Mathematical Surveys and Monographs,
233, Amer.
Math. Soc., Providence, RI, 2018.


\bibitem{Kr_20} N.V. Krylov,  {\em
Linear and fully nonlinear elliptic equations with $L_{d}$-drift\/},
  Comm. PDE, Vol. 45 (2020), No.  12, 1778--1798.


 \bibitem{Kr_21b}  N.V. Krylov,  {\em
Elliptic equations with VMO a,
  b$\,\in L_{d}$, and c$\,\in L_{d/2}$\/}, Trans. Amer. Math.
Sci., Vol. 374  (2021), No. 4, 2805-2822.

\bibitem{Kr_21}  N.V. Krylov,
{\em On stochastic equations with drift in $L_{d}$\/}, arXiv:2001.04008.

\bibitem{Kr_21a}  N.V. Krylov,
{\em On   diffusion
 processes with drift in a Morrey class containing  $L_{d+2}$\/}, arXiv:2104.05603

 \bibitem{St_70} E. Stein, ``Singular integrals and
differentiability properties of functions'',
Princeton University Press, Princeton, NJ, 1970.



\bibitem{Wi_09} N. Winter, {\em $W^{2,p}$ and $W^{1,p}$-estimates
at the boundary for solutions of fully nonlinear,
uniformly elliptic equations\/},
 Z. Anal. Anwend., Vol. 28  (2009), No. 2, 129--164.


\end{thebibliography}
\end{document}